\documentclass[12pt,oneside,english]{amsart}
\usepackage[T1]{fontenc}
\usepackage[latin9]{inputenc}
\usepackage{geometry}
\usepackage{amsthm}
\usepackage{amssymb}
\usepackage{graphicx}
\usepackage{setspace}
\makeatletter
\numberwithin{equation}{section}
\numberwithin{figure}{section}
\theoremstyle{plain}
\newtheorem{thm}{\protect\theoremname}
\theoremstyle{definition}
\newtheorem{example}[thm]{\protect\examplename}

\makeatother

\usepackage{babel}
\providecommand{\examplename}{Example}
\providecommand{\theoremname}{Theorem}

\begin{document}

\title{Multivariate Force of Mortality}

\maketitle
(Appeared in the journal \emph{Demography India} in 2015)

\begin{center}

\textbf{\large{}Swagata Mitra}{\large\par}

Data Processing Division 

National Sample Survey Office, Kolkata 

Ministry of Statistics \& Programme Implementation

Government of India, Kolkata INDIA 700108

Email: swagatamitrastat09@gmail.com

$ $

\textbf{\large{}Pratyush Singh}{\large\par}

Indian Statistical Institute, 

Kolkata, INDIA 700108

Email: pratyush.isi@gmail.com

$ $

\textbf{\large{}Arni S.R. Srinivasa Rao}\footnote{{\large{}Corresponding author}}\textbf{\large{},}\footnote{{\large{}Funding for this work was received by ASRSR (Principal Investigator)
when he was a permanent faculty at Indian Statistical Institute (ISI),
Kolkata for the Project Entitled \textquotedblleft Multiple Decrement
Tables in Population Health Insurance policies: Deterministic Approaches\textquotedblright ,
funded by the ISI for the period 2010-2012. SM worked with ASRSR under
the project as a project-linked person and PS worked with ASRSR under
the project as a M.Stat student of ISI during the period 2010-2012.}\textbf{\large{} }}{\large\par}

Department of Biostatistics and Epidemiology,

Medical College of Georgia and College of Science Mathematics,

Augusta University, 1120 15th Street, Augusta, GA 30912, USA

Email: arrao@augusta.edu

\end{center}
\begin{abstract}
In usual demographic analysis, force of mortality is a function of
one variable, that is, of age. In this article bi-variate and multivariate
force of mortality functions are introduced for the first time to
explain mortality differentials. The pattern of mortality in a population
is one of the strong influencing factors in determining the life expectancies
at various ages in the population. Considering univariate functions
of age only to understand the human mortality data without associating
with other variables could lead to incomplete analysis. The reasons
behind declining forces of mortality globally could be studied using
the proposed functions. Other applications of multivariate forces
of mortality could be in actuarial sciences. 
\end{abstract}

\keywords{Keywords: Hazard rate, PDEs, Taylor's expansions. }

\subjclass[2000]{MSC:91D20, 65.0X}

\tableofcontents{}

\section{Introduction}

One variable force of mortality, $\mu(x)$, with respect to age, $x$
of an individual is one of the central topics of study in the actuarial
mathematics and there are very useful discussions are available on
univariate or one variable force of mortality (see for example, \cite{SmithAMM1948,TurnerJRSSA2010,Finkelstein2003}).
It is often termed as instantaneous rate of death at an age $x.$
Suppose, forces of mortality is measured on two variables $(x,y)$,
(one being age, other can be some influencing factor on mortality),
then if we plot the force of mortality $\mu(x,y)$ on the $xy-$plane,
then for an arbitrary point $(x_{0},y_{0})$, we can write,

\begin{eqnarray}
\mu\left(x_{0},y_{0}\right) & = & \frac{1}{\rho(\Omega)}\int\int_{\Omega}\mu(x,y)d\rho\label{eq:0-0}
\end{eqnarray}

where $\Omega$ is a region such that $\left(x_{0},y_{0}\right)\in\Omega$
and $\rho(\Omega)$ is area of the region. Since $\mu$ is continuous
in the univariate case, if we assume the same holds for $\mu\left(x,y\right)$
in $\Omega$, it has an upper bound say, $\mu_{1}$ and a lower bound
say, $\mu_{0}$ in the region such that 

\begin{eqnarray*}
\mu_{0}\leq\frac{1}{\rho(\Omega)}\int\int_{\Omega}\mu(x,y)d\rho\leq\mu_{1}.
\end{eqnarray*}

Multivariate force of mortality functions can help in better understanding
of future longevity and causes of decline in mortality rates. There
are studies which consider mortality decline or longevity projections
of humans with respect to age only (for example, see \cite{TurnerJRSSA2010,Finkelstein2003}).
Such studies can be handled using univariate analysis of standard
force of mortality functions. However such studies can be extended
to incorporate several variables that can explain decline in mortality
rates or causes of increase in longevity using multivariate forces
of mortality functions. Further analytical properties such as rate
of change in forces of mortalities of bi-variate force of mortality
functions are found in the next sections. We will begin with two examples
and we illustrate their numerical properties.
\begin{example}
Suppose $l(x,y)=1-x^{a}y^{b/\sqrt{K}}/K$ for $a,b,K\in\mathbb{Z}^{+}.$
See the definitions of $\mu_{x}(x,y)$

and $\mu_{y}(x,y)$ in eq. (\ref{eq:mu-5}) and eq. (\ref{eq:mu-6})
in the section 3. We have,

\begin{eqnarray*}
\mu_{x}(x,y) & = & \frac{ax^{a-1}y^{b/\sqrt{K}}}{K-x^{a}y^{b/\sqrt{K}}}\\
\mu_{y}(x,y) & = & \frac{bx^{a}y^{(b/\sqrt{K})-1}}{\sqrt{K}\left(K-x^{a}y^{b/\sqrt{K}}\right)}
\end{eqnarray*}
\end{example}

$ $
\begin{example}
Suppose $l(x,y)=a^{\sqrt{y}}b^{x^{3}}$ for $10\leq x\leq80,$ $1\leq y\leq5,$
$a=5,$ $b=7.$ $\mu_{x}(x,y)=-3x^{2}\log(b)$ and $\mu_{y}(x,y)=-\log(a)/2\sqrt{y}.$
\end{example}

\begin{figure}
\includegraphics[scale=0.8]{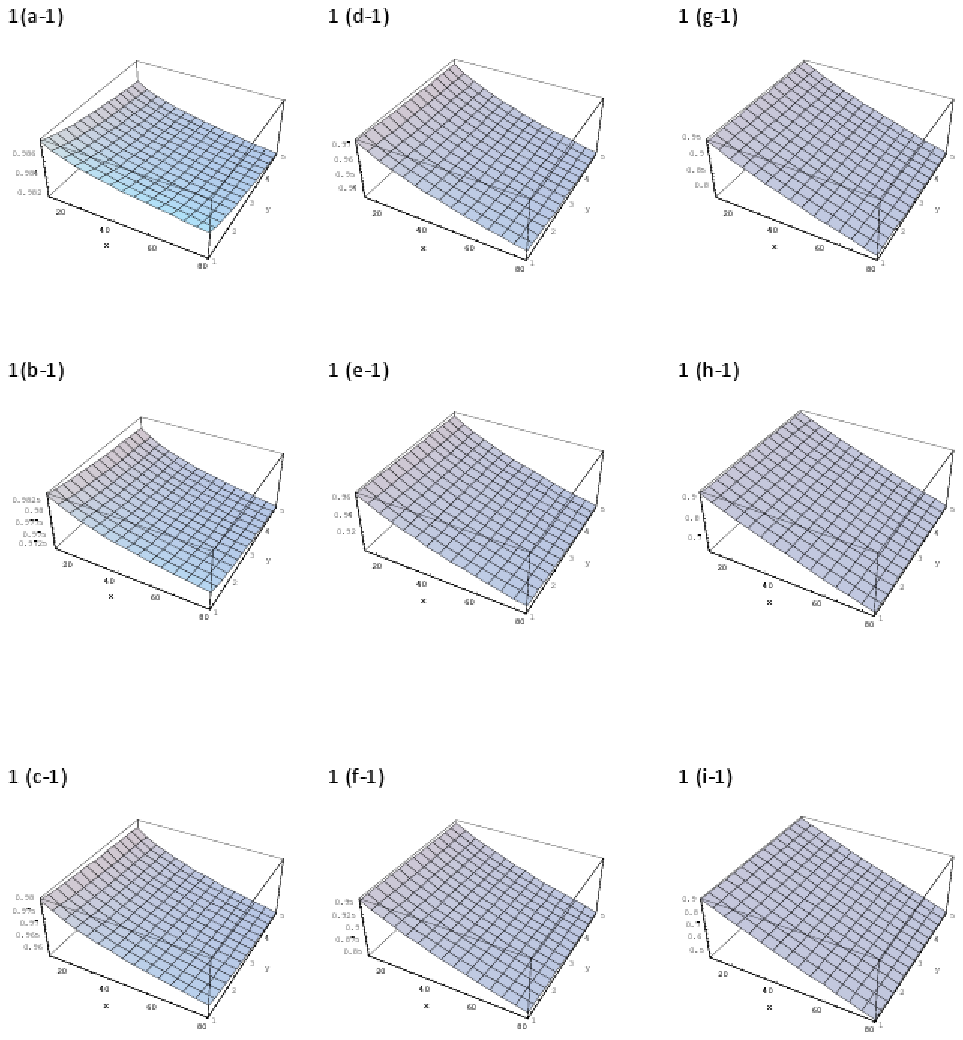}

\caption{$l(x,y)$ in Example 1. For various combinations of values of \textbf{$a$}
and $b$ we have drawn $1(a-1)$ to $1(i-1)$ by fixing $K=100,$
$10\protect\leq x\protect\leq80$ and $1\protect\leq y\protect\leq5.$
Following are the combinations of $a$ and $b$ for each figure: $1(a-1):a=0.1,b=0.9$,
$1(b-1):a=0.2,b=0.8,$ $1(c-1):a=0.3,b=0.7,$ $1(d-1):a=0.4,b=0.6,$
$1(e-1):a=0.5,b=0.5,$ $1(f-1):a=0.6,b=0.4,$ $1(g-1):a=0.7,b=0.3,$
$1(h-1):a=0.8,b=0.2,$ $1(i-1):a=0.9,b=0.1.$ }

\end{figure}

\begin{figure}
\includegraphics[scale=0.8]{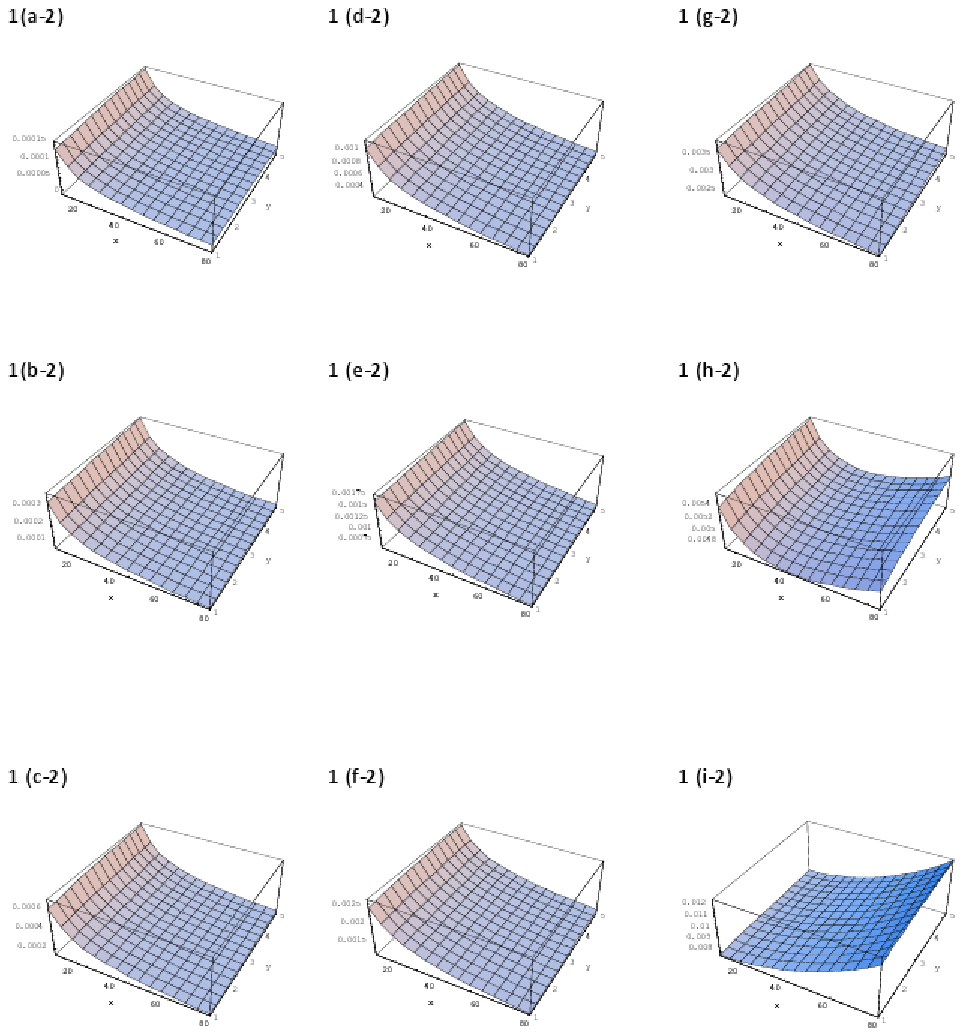}

\caption{$\mu_{x}(x,y)$ in Example 1. For various combinations of values of
\textbf{$a$} and $b$ we have drawn $1(a-1)$ to $1(i-1)$ by fixing
$K=100,$ $10\protect\leq x\protect\leq80$ and $1\protect\leq y\protect\leq5.$
Following are the combinations of $a$ and $b$ for each figure: $1(a-1):a=0.1,b=0.9$,
$1(b-1):a=0.2,b=0.8,$ $1(c-1):a=0.3,b=0.7,$ $1(d-1):a=0.4,b=0.6,$
$1(e-1):a=0.5,b=0.5,$ $1(f-1):a=0.6,b=0.4,$ $1(g-1):a=0.7,b=0.3,$
$1(h-1):a=0.8,b=0.2,$ $1(i-1):a=0.9,b=0.1.$ }

\end{figure}

\begin{figure}
\includegraphics[scale=0.8]{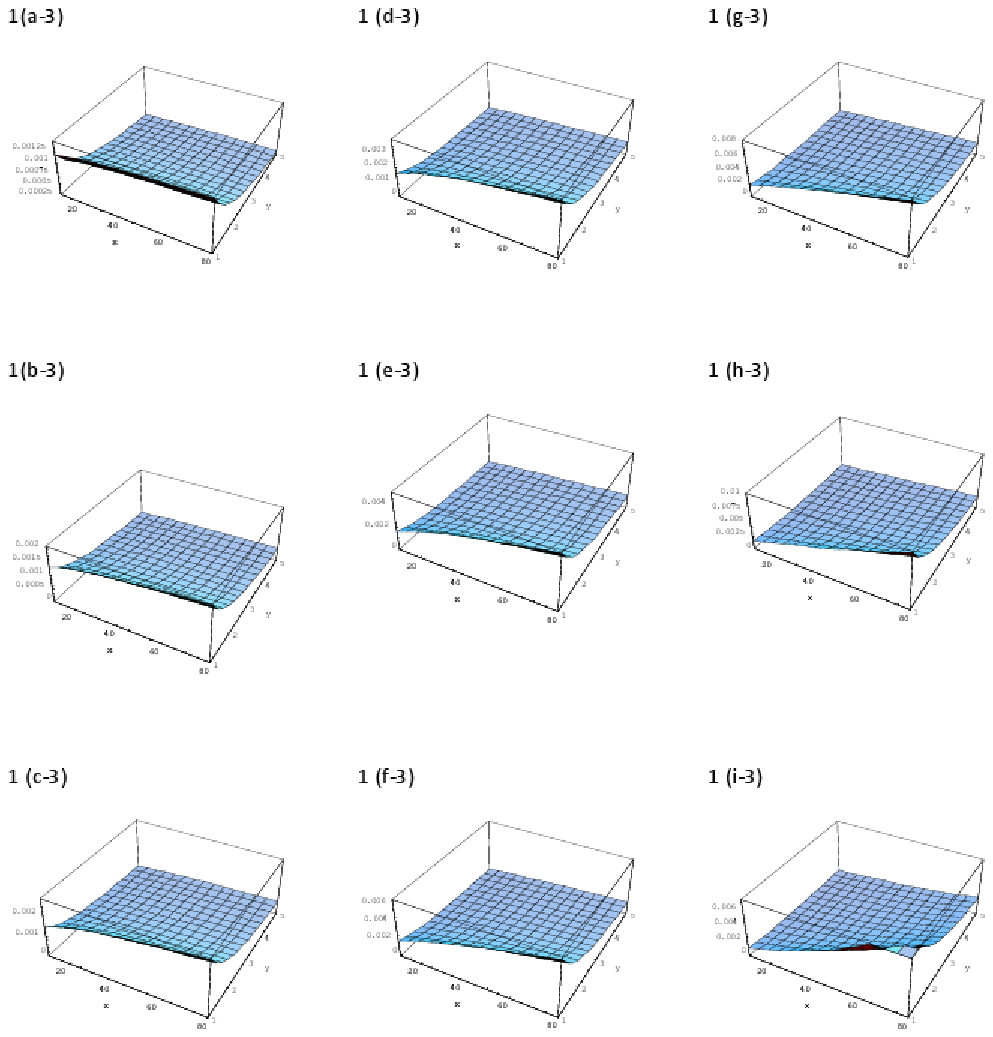}

\caption{$\mu_{y}(x,y)$ in Example 1. For various combinations of values of
\textbf{$a$} and $b$ we have drawn $1(a-1)$ to $1(i-1)$ by fixing
$K=100,$ $10\protect\leq x\protect\leq80$ and $1\protect\leq y\protect\leq5.$
Following are the combinations of $a$ and $b$ for each figure: $1(a-1):a=0.1,b=0.9$,
$1(b-1):a=0.2,b=0.8,$ $1(c-1):a=0.3,b=0.7,$ $1(d-1):a=0.4,b=0.6,$
$1(e-1):a=0.5,b=0.5,$ $1(f-1):a=0.6,b=0.4,$ $1(g-1):a=0.7,b=0.3,$
$1(h-1):a=0.8,b=0.2,$ $1(i-1):a=0.9,b=0.1.$ }

\end{figure}

Force of mortality functions also associated with continuous life
table functions \cite{SmithAMM1948}. Life table is a mathematical
model describing how individuals born at same time (for example, a
cohort of new born babies) survive over the years at various ages
until the last individual dies. Life table is constructed based on
present or past mortality pattern (i.e. mortality rates at each age,
(say $x$) per fixed number of individuals in the same age $x$ per
year or for the year $(0,t)$) in the population and assumed that
this pattern will remain the same until the individual at last age
dies. This table can be used to construct synthetic population at
each age $x$ at time $t$ or for time interval $(0,t)$, it do not
have mechanism to take care of future changes in the mortality pattern
after $t.$ Suppose $l(x)$ denote the number of individuals at age
$x$ out of $l(0)$ newly born individuals in a life table with continuous
partial derivatives up to order $(k+1)$. Then the number of individuals
at ages, $x+\Delta x$ and $x-\Delta x$ are denoted by $l(x+\Delta x)$
and $l(x-\Delta x)$ can be obtained from Taylor series expansion
evaluating at $x_{0}$ as 

\begin{eqnarray}
l(x_{0}+\Delta x) & = & l(x_{0})+\Delta xl'(x_{0})+(\Delta x)^{2}\frac{l^{(2)}(x_{0})}{2!}+(\Delta x)^{3}\frac{l^{(3)}(x_{0})}{3!}+...\nonumber \\
 &  & \qquad+(\Delta x)^{k}\frac{l^{(k)}(x_{0})}{k!}+\int_{x_{0}}^{x}\frac{(\Delta x)^{k}l^{(k+1)}(t)}{k!}dt\label{eq:1.1}\\
l(x_{0}-\Delta x) & = & l(x_{0})-\Delta xl'(x_{0})+(\Delta x)^{2}\frac{l^{(2)}(x_{0})}{2!}-(\Delta x)^{3}\frac{l^{(3)}(x_{0})}{3!}+...\nonumber \\
 &  & \qquad+(-1)^{k}(\Delta x)^{k}\frac{l^{(k)}(x_{0})}{k!}+\int_{x_{0}}^{x}\frac{(\Delta x)^{k}l^{(k+1)}(t)}{k!}dt\label{eq:1.2}
\end{eqnarray}

Here $x$ is the fixed age between 0 and $\omega$, the maximum age
of life. The number of survivors at age $x$ with some other relevant
factor (for example, marital status, education level, climate, food
habits, geographic region etc) $y$ evaluated at $\left(x_{0},y_{0}\right)$
can be obtained form two variable Taylor expansion. Assuming continuous
partial derivatives for $l\left(x_{0}+\Delta x,y_{0}+\Delta y\right)$
up to order $3$, we can the expansion of two variable survival functions
as follows:

\begin{eqnarray}
l\left(x_{0}+\Delta x,y_{0}+\Delta y\right) & = & l\left(x_{0},y_{0}\right)+\Delta x\left\{ \frac{\partial l}{\partial x}\left(x_{0},y_{0}\right)\right\} +\Delta y\left\{ \frac{\partial l}{\partial y}\left(x_{0},y_{0}\right)\right\} +\frac{(\Delta x)^{2}}{2!}\left\{ \frac{\partial^{2}l}{\partial x^{2}}\left(x_{0},y_{0}\right)\right\} \nonumber \\
\nonumber \\
 &  & +\frac{(\Delta y)^{2}}{2!}\left\{ \frac{\partial^{2}l}{\partial y^{2}}\left(x_{0},y_{0}\right)\right\} +\Delta x\Delta y\left\{ \frac{\partial^{2}l}{\partial x\partial y}\left(x_{0},y_{0}\right)\right\} +\frac{(\Delta x)^{3}}{3!}\left\{ \frac{\partial^{3}l}{\partial x^{3}}\left(x_{0},y_{0}\right)\right\} \nonumber \\
\nonumber \\
 &  & +\frac{(\Delta y)^{3}}{3!}\left\{ \frac{\partial^{3}l}{\partial y^{3}}\left(x_{0},y_{0}\right)\right\} +\frac{\left(\Delta x\right)^{2}\Delta y}{2}\left\{ \frac{\partial^{3}l}{\partial x^{2}\partial y}\left(x_{0},y_{0}\right)\right\} \nonumber \\
 &  & +\frac{\Delta x\left(\Delta y\right)^{2}}{2}\left\{ \frac{\partial^{3}l}{\partial x\partial y^{2}}\left(x_{0},y_{0}\right)\right\} +\cdots+\cdots\label{eq:1-3}
\end{eqnarray}

\begin{eqnarray}
l\left(x_{0}-\Delta x,y_{0}-\Delta y\right) & = & l\left(x_{0},y_{0}\right)-\Delta x\left\{ \frac{\partial l}{\partial x}\left(x_{0},y_{0}\right)\right\} -\Delta y\left\{ \frac{\partial l}{\partial y}\left(x_{0},y_{0}\right)\right\} +\frac{(\Delta x)^{2}}{2!}\left\{ \frac{\partial^{2}l}{\partial x^{2}}\left(x_{0},y_{0}\right)\right\} \nonumber \\
\nonumber \\
 &  & +\frac{(\Delta y)^{2}}{2!}\left\{ \frac{\partial^{2}l}{\partial y^{2}}\left(x_{0},y_{0}\right)\right\} +\Delta x\Delta y\left\{ \frac{\partial^{2}l}{\partial x\partial y}\left(x_{0},y_{0}\right)\right\} -\frac{(\Delta x)^{3}}{3!}\left\{ \frac{\partial^{3}l}{\partial x^{3}}\left(x_{0},y_{0}\right)\right\} \nonumber \\
\nonumber \\
 &  & -\frac{(\Delta y)^{3}}{3!}\left\{ \frac{\partial^{3}l}{\partial y^{3}}\left(x_{0},y_{0}\right)\right\} -\frac{\left(\Delta x\right)^{2}\Delta y}{2}\left\{ \frac{\partial^{3}l}{\partial x^{2}\partial y}\left(x_{0},y_{0}\right)\right\} \nonumber \\
 &  & -\frac{\Delta x\left(\Delta y\right)^{2}}{2}\left\{ \frac{\partial^{3}l}{\partial x\partial y^{2}}\left(x_{0},y_{0}\right)\right\} +\cdots+\cdots\label{eq:1-4}
\end{eqnarray}

\section{Analysis of first order equations}

We try to analyze univariate force of mortality functions by considering
higher order derivatives are continuous. Assuming $x\rightarrow x_{0}$,
in eq. (\ref{eq:1.1}) and eq. (\ref{eq:1.2}), we will have $\int_{x_{0}}^{x}\left((\Delta x)^{k}l^{(k+1)}(t)/k!\right)dt\rightarrow0.$
Suppose $f^{(n)}(x_{0}+\Delta x)\mbox{ and }f^{(n)}(x_{0}-\Delta x)$
denote Taylor expansion equations when ignoring the $(n+1)^{th}$
derivatives and beyond for $n=2,3,...$ in eq. (\ref{eq:1.1}) and
eq. (\ref{eq:1.2}), then by sequentially ignoring the terms beginning
from the term $\frac{(\Delta x)^{n}}{n!}l^{(n)}(x_{0})$ in eq. (\ref{eq:1.1})
and eq. (\ref{eq:1.2}) for $n=2,3,...$ , we will obtain following
equations:

\begin{eqnarray}
f^{(1)}(x_{0}+\Delta x)-f^{(1)}(x_{0}-\Delta x) & = & 2\Delta xl'(x_{0})\label{eq:l1}\\
f^{(3)}(x_{0}+\Delta x)-f^{(3)}(x_{0}-\Delta x) & = & 2\Delta xl'(x_{0})+\frac{2(\Delta x)^{3}}{3!}l^{(3)}(x_{0})\label{eq:l3}
\end{eqnarray}

\begin{eqnarray}
f^{(5)}(x_{0}+\Delta x)-f^{(5)}(x_{0}-\Delta x) & = & 2\Delta xl'(x_{0})+\frac{2(\Delta x)^{3}}{3!}l^{(3)}(x_{0})+\frac{2(\Delta x)^{5}}{5!}l^{(5)}(x_{0})\nonumber \\
\vdots &  & \vdots\nonumber \\
f^{(2k-1)}(x_{0}+\Delta x)-f^{(2k-1)}(x_{0}-\Delta x) & = & 2\Delta xl'(x_{0})+\frac{2(\Delta x)^{3}}{3!}l^{(3)}(x_{0})+...+\frac{2(\Delta x)^{(2k-1)}}{(2k-1)!}l^{(2k-1)}(x_{0})\nonumber \\
\vdots &  & \vdots\label{eq:2}
\end{eqnarray}

Therefore,

$\sum_{j=1}^{k}\left[\left\{ f^{(2j+1)}(x_{0}+\Delta x)-f^{(2j+1)}(x_{0}-\Delta x)\right\} -\left\{ f^{(2j-1)}(x_{0}+\Delta x)-f^{(2j-1)}(x_{0}-\Delta x)\right\} \right]=$

\begin{eqnarray}
 &  & \frac{2(\Delta x)^{3}}{3!}l^{(3)}(x_{0})+\frac{2(\Delta x)^{5}}{5!}l^{(5)}(x_{0})+...+\frac{2(\Delta x)^{(2k+1)}}{(2k+1)!}l^{(2k+1)}(x_{0})\label{eq:3}
\end{eqnarray}

Hence we will obtain,

$\sum_{j=1}^{k}\left[\left\{ f^{(2j+1)}(x_{0}+\Delta x)-f^{(2j+1)}(x_{0}-\Delta x)\right\} -\left\{ f^{(2j-1)}(x_{0}+\Delta x)-f^{(2j-1)}(x_{0}-\Delta x)\right\} \right]=$

\begin{eqnarray}
 &  & \left\{ f^{(2k+1)}(x_{0}+\Delta x)-f^{(2k+1)}(x_{0}-\Delta x)\right\} -\left\{ f^{(1)}(x_{0}+\Delta x)-f^{(1)}(x_{0}-\Delta x)\right\} \label{eq:4}
\end{eqnarray}

Using the relation $l^{(1)}(x_{0})=\frac{1}{2\Delta x}\left[f^{(1)}(x_{0}+\Delta x)-f^{(1)}(x_{0}-\Delta x)\right]$,
we write,

\begin{eqnarray}
\left.\frac{dl(x)}{dx}\right|_{x=x_{0}} & = & \frac{1}{2\Delta x}\Biggl[\left\{ f^{(2k+1)}(x_{0}+\Delta x)-f^{(2k+1)}(x_{0}-\Delta x)\right\} -\nonumber \\
 &  & \left.\sum_{j=1}^{k}\left[\left\{ f^{(2j+1)}(x_{0}+\Delta x)-f^{(2j+1)}(x_{0}-\Delta x)\right\} -\left\{ f^{(2j-1)}(x_{0}+\Delta x)-f^{(2j-1)}(x_{0}-\Delta x)\right\} \right]\right]\label{eq:5}
\end{eqnarray}

Dividing the eq. (\ref{eq:l1}) by the term $2\Delta x$ on both the
sides and integrating it from age $x$ to $x+m$, we obtain,

\begin{eqnarray}
\int_{x}^{x+m}\frac{\left\{ f^{(1)}(y_{0}+\Delta y)-f^{(1)}(y_{0}-\Delta y)\right\} }{2\Delta y}dy & = & \int_{x}^{x+m}\left.\frac{d}{dy}l(y)\right|_{y=y_{0}}dy\nonumber \\
 & = & -\int_{x}^{x+m}\left.-\frac{1}{l(y)}\frac{d}{dy}l(y)\right|_{y=y_{0}}\left.l(y)\right|_{y=y_{0}}dy\nonumber \\
 & = & -\int_{x}^{x+m}\left.\mu(y)\right|_{y=y_{0}}\left.l(y)\right|_{y=y_{0}}dy\label{eq:6}
\end{eqnarray}

where, $\mu(x)$, the force of mortality function, which is defined
as $\left(-1/l(x)\right)(d/dx)l(x).$ 

Dividing the eq. (\ref{eq:l3}) by the term $2\Delta x$ on both the
sides and integrating it from age $x$ to $x+m$, we obtain,

\begin{eqnarray}
\int_{x}^{x+m}\left[\frac{\left\{ f^{(3)}(y_{0}+\Delta y)-f^{(3)}(y_{0}-\Delta y)\right\} }{2\Delta y}-\frac{(\Delta y)^{2}}{3!}l^{(3)}(y{}_{0})\right]dy & = & -\int_{x}^{x+m}\left.\mu(y)\right|_{y=y_{0}}\left.l(y)\right|_{y=y_{0}}dy\label{eq:7}
\end{eqnarray}

Now, multiplying $\left.l(x)\right|_{x=x_{0}}$ and $\left.(-1/l(x))\right|_{x=x_{0}}$
to the eq. (\ref{eq:5}), and integrating between ages $x$ and $x+m$,
we will obtain,

\begin{eqnarray*}
-\int_{x}^{x+m}\left.-\frac{1}{l(y)}\frac{d}{dy}l(y)\right|_{y=y_{0}}\left.l(y)\right|_{y=y_{0}}dy & = & \int_{x}^{x+m}-\frac{1}{2\Delta y}\Biggl[\left\{ f^{(2k+1)}(y_{0}+\Delta y)-f^{(2k+1)}(y_{0}-\Delta y)\right\} \\
 &  & -\sum_{j=1}^{k}\left[\left\{ f^{(2j+1)}(y_{0}+\Delta y)-f^{(2j+1)}(y_{0}-\Delta y)\right\} \right.\\
 &  & \left.\left.-\left\{ f^{(2j-1)}(y_{0}+\Delta y)-f^{(2j-1)}(y_{0}-\Delta y)\right\} \right]\right]
\end{eqnarray*}

\begin{eqnarray}
 & = & -\int_{x}^{x+m}\left.\frac{d}{dy}l(y)\right|_{y=y_{0}}dy\nonumber \\
 & = & ml'(y_{0})\label{eq:9}
\end{eqnarray}

If there are any deaths during the age $(x_{0}+\Delta x)$ to $(x_{0}-\Delta x)$,
then $f^{(1)}(x_{0}+\Delta x)<f^{(1)}(x_{0}-\Delta x).$ In the absence
of deaths, we have $f^{(1)}(x_{0}+\Delta x)=f^{(1)}(x_{0}-\Delta x)$
and $\left.\frac{dl(x)}{dx}\right|_{x=x_{0}}<0.$ This argument is
true for difference of other higher order expressions. 

\section{Analysis of second order equations}

Suppose $g^{(n)}(x_{0}+\Delta x,y_{0}+\Delta x)\mbox{ and }g^{(n)}(x_{0}-\Delta x,y_{0}-\Delta y)$
denote equations when ignoring the terms from $(n+1)^{th}$ partial
derivatives and beyond for $n=2,3,...$ in eq. (\ref{eq:1-3}) and
eq. (\ref{eq:1-4}), then by sequentially ignoring the terms beginning
from the $n^{th}$ order partial derivative terms in eq. (\ref{eq:1-3})
and eq. (\ref{eq:1-4}) for $n=2,3,...$ , we will obtain following
equations:

\begin{eqnarray*}
g^{(1)}\left(x_{0}+\Delta x,y_{0}+\Delta y\right)-g^{(1)}\left(x_{0}-\Delta x,y_{0}-\Delta y\right) & = & 2\Delta x\left\{ \frac{\partial l}{\partial x}\left(x_{0},y_{0}\right)\right\} +2\Delta y\left\{ \frac{\partial l}{\partial y}\left(x_{0},y_{0}\right)\right\} \\
g^{(3)}\left(x_{0}+\Delta x,y_{0}+\Delta y\right)-g^{(3)}\left(x_{0}-\Delta x,y_{0}-\Delta y\right) & = & 2\Delta x\left\{ \frac{\partial l}{\partial x}\left(x_{0},y_{0}\right)\right\} +2\Delta y\left\{ \frac{\partial l}{\partial y}\left(x_{0},y_{0}\right)\right\} \\
 &  & +\frac{(\Delta x)^{3}}{3!}\left\{ \frac{\partial^{3}l}{\partial x^{3}}\left(x_{0},y_{0}\right)\right\} +\frac{(\Delta y)^{3}}{3!}\left\{ \frac{\partial^{3}l}{\partial y^{3}}\left(x_{0},y_{0}\right)\right\} 
\end{eqnarray*}

\begin{eqnarray*}
 &  & +\frac{\left(\Delta x\right)^{2}\Delta y}{2}\left\{ \frac{\partial^{3}l}{\partial x^{2}\partial y}\left(x_{0},y_{0}\right)\right\} \\
 &  & +\frac{\Delta x\left(\Delta y\right)^{2}}{2}\left\{ \frac{\partial^{3}l}{\partial x\partial y^{2}}\left(x_{0},y_{0}\right)\right\} \\
\vdots &  & \vdots
\end{eqnarray*}

\begin{eqnarray*}
g^{(2k-1)}\left(x_{0}+\Delta x,y_{0}+\Delta y\right)-g^{(2k-1)}\left(x_{0}-\Delta x,y_{0}-\Delta y\right) & = & 2\Delta x\left\{ \frac{\partial l}{\partial x}\left(x_{0},y_{0}\right)\right\} +2\Delta y\left\{ \frac{\partial l}{\partial y}\left(x_{0},y_{0}\right)\right\} +...\\
 &  & +\frac{2\left(\Delta x\right)^{2k-j-1}\left(\Delta y\right)^{j}}{\left(2k-1\right)!}\left[\sum_{j=0}^{2k-1}\left(\begin{array}{c}
2k-1\\
j
\end{array}\right)\times\right.\\
 &  & \left.\left\{ \frac{\partial^{(2k-1)}l}{\partial x^{(2k-j-1)}\partial y^{j}}\left(x_{0},y_{0}\right)\right\} \right]\\
\vdots &  & \vdots
\end{eqnarray*}

Therefore,

\begin{eqnarray*}
\sum_{j=1}^{k}\left[\left\{ g^{(2k+1)}\left(x_{0}+\Delta x,y_{0}+\Delta y\right)-g^{(2k+1)}\left(x_{0}-\Delta x,y_{0}-\Delta y\right)\right\} \right.\\
-\left.\left\{ g^{(2k-1)}\left(x_{0}+\Delta x,y_{0}+\Delta y\right)-g^{(2k-1)}\left(x_{0}-\Delta x,y_{0}-\Delta y\right)\right\} \right] & = & \frac{(\Delta x)^{3}}{3!}\left\{ \frac{\partial^{3}l}{\partial x^{3}}\left(x_{0},y_{0}\right)\right\} 
\end{eqnarray*}

\begin{eqnarray}
 &  & +\frac{(\Delta y)^{3}}{3!}\left\{ \frac{\partial^{3}l}{\partial y^{3}}\left(x_{0},y_{0}\right)\right\} +\frac{\left(\Delta x\right)^{2}\Delta y}{2}\left\{ \frac{\partial^{3}l}{\partial x^{2}\partial y}\left(x_{0},y_{0}\right)\right\} +\frac{\Delta x\left(\Delta y\right)^{2}}{2}\left\{ \frac{\partial^{3}l}{\partial x\partial y^{2}}\left(x_{0},y_{0}\right)\right\} +\nonumber \\
 &  & \cdots+\frac{2\left(\Delta x\right)^{2k-j+1}\left(\Delta y\right)^{j}}{\left(2k+1\right)!}\left[\sum_{j=0}^{2k+1}\left(\begin{array}{c}
2k+1\\
j
\end{array}\right)\left\{ \frac{\partial^{(2k+1)}l}{\partial x^{(2k-j+1)}\partial y^{j}}\left(x_{0},y_{0}\right)\right\} \right]\label{eq:master-2}
\end{eqnarray}

Force of mortality for the two variables $\left(x,y\right)$ is evaluated
at the point $\left(x_{0},y_{0}\right)$ using partial derivatives
as follows:

\begin{eqnarray}
\frac{\partial\mu\left(x,y\right)}{\partial x} & = & \lim_{\Delta x\rightarrow0}\frac{\mu\left(x_{0}+\Delta x,y_{0}\right)-\mu\left(x_{0,}y_{0}\right)}{\Delta x}\label{eq:mu-1}\\
\nonumber \\
\frac{\partial\mu\left(x,y\right)}{\partial y} & = & \lim_{\Delta y\rightarrow0}\frac{\mu\left(x_{0},y_{0}+\Delta y\right)-\mu\left(x_{0,}y_{0}\right)}{\Delta y}\label{eq:mu-2}
\end{eqnarray}

where, we define, 
\begin{eqnarray}
\mu\left(x_{0}+\Delta x,y_{0}\right) & = & -\frac{1}{l\left(x_{0}+\Delta x,y_{0}\right)}\frac{\partial l\left(x_{0}+\Delta x,y_{0}\right)}{\partial x}\label{eq:mu-3}\\
\mu\left(x_{0},y_{0}+\Delta y\right) & = & -\frac{1}{l\left(x_{0},y_{0}+\Delta y\right)}\frac{\partial l\left(x_{0},y_{0}+\Delta y\right)}{\partial y}\label{eq:mu-4}\\
\mu_{x}\left(x_{0,}y_{0}\right) & = & -\frac{1}{l\left(x_{0},y_{0}\right)}\frac{\partial l\left(x_{0},y_{0}\right)}{\partial x}\mbox{ for eq. (\ref{eq:mu-1})}\label{eq:mu-5}\\
\mu_{y}\left(x_{0,}y_{0}\right) & = & -\frac{1}{l\left(x_{0},y_{0}\right)}\frac{\partial l\left(x_{0},y_{0}\right)}{\partial y}\mbox{ for eq. (\ref{eq:mu-2})}\label{eq:mu-6}
\end{eqnarray}

Using eq. (\ref{eq:1-3}) and as $\Delta y\rightarrow0$, we obtain,

\begin{eqnarray}
l\left(x_{0}+\Delta x,y_{0}\right) & = & l\left(x_{0},y_{0}\right)+\Delta x\left\{ \frac{\partial l}{\partial x}\left(x_{0},y_{0}\right)\right\} +\frac{(\Delta x)^{2}}{2!}\left\{ \frac{\partial^{2}l}{\partial x^{2}}\left(x_{0},y_{0}\right)\right\} +\frac{(\Delta x)^{3}}{3!}\left\{ \frac{\partial^{3}l}{\partial x^{3}}\left(x_{0},y_{0}\right)\right\} \nonumber \\
\label{eq:3.8,1 for mu-3}
\end{eqnarray}

Using eq. (\ref{eq:1-4}) and as $\Delta x\rightarrow0$, we obtain,

\begin{eqnarray}
l\left(x_{0},y_{0}+\Delta y\right) & = & l\left(x_{0},y_{0}\right)+\Delta y\left\{ \frac{\partial l}{\partial y}\left(x_{0},y_{0}\right)\right\} +\frac{(\Delta y)^{2}}{2!}\left\{ \frac{\partial^{2}l}{\partial y^{2}}\left(x_{0},y_{0}\right)\right\} +\frac{(\Delta y)^{3}}{3!}\left\{ \frac{\partial^{3}l\left(x_{0},y_{0}\right)}{\partial y^{3}}\right\} \nonumber \\
\label{eq:3.9,1 for mu-4}
\end{eqnarray}

Therefore,

\begin{eqnarray}
\frac{\partial l}{\partial x}\left(x_{0}+\Delta x,y_{0}\right) & = & \left\{ \frac{\partial l}{\partial x}\left(x_{0},y_{0}\right)\right\} +\Delta x\left\{ \frac{\partial^{2}l}{\partial x^{2}}\left(x_{0},y_{0}\right)\right\} +\frac{(\Delta x)^{2}}{2!}\left\{ \frac{\partial^{3}l}{\partial x^{3}}\left(x_{0},y_{0}\right)\right\} +\frac{(\Delta x)^{3}}{3!}\left\{ \frac{\partial^{4}l}{\partial x^{4}}\left(x_{0},y_{0}\right)\right\} \nonumber \\
\label{eq:3.10,do for mu-3}\\
\frac{\partial l}{\partial y}\left(x_{0},y_{0}+\Delta y\right) & = & \left\{ \frac{\partial l}{\partial y}\left(x_{0},y_{0}\right)\right\} +\Delta y\left\{ \frac{\partial^{2}l}{\partial y^{2}}\left(x_{0},y_{0}\right)\right\} +\frac{(\Delta y)^{2}}{2!}\left\{ \frac{\partial^{3}l}{\partial y^{3}}\left(x_{0},y_{0}\right)\right\} +\frac{(\Delta y)^{3}}{3!}\left\{ \frac{\partial^{4}l}{\partial y^{4}}\left(x_{0},y_{0}\right)\right\} \nonumber \\
\label{eq:3.11,do for mu-4}
\end{eqnarray}

Let us now derive the equation of the type eq. (\ref{eq:master-2})
with the conditions $\Delta y\rightarrow0$ and $\Delta x\rightarrow0$
and extending these equations up to the general term. Suppose, $\Delta y\rightarrow0$
in the equations (\ref{eq:1-3}) and (\ref{eq:1-4}) , then 

\begin{eqnarray}
l\left(x_{0}+\Delta x,y_{0}\right) & = & l\left(x_{0},y_{0}\right)+\Delta x\left\{ \frac{\partial l}{\partial x}\left(x_{0},y_{0}\right)\right\} +\frac{(\Delta x)^{2}}{2!}\left\{ \frac{\partial^{2}l}{\partial x^{2}}\left(x_{0},y_{0}\right)\right\} +\cdots\nonumber \\
 &  & \qquad+\frac{(\Delta x)^{n}}{n!}\left\{ \frac{\partial^{n}l}{\partial x^{n}}\left(x_{0},y_{0}\right)\right\} +\cdots\label{eq:2-1}\\
l\left(x_{0}-\Delta x,y_{0}\right) & = & l\left(x_{0},y_{0}\right)-\Delta x\left\{ \frac{\partial l}{\partial x}\left(x_{0},y_{0}\right)\right\} +\frac{(\Delta x)^{2}}{2!}\left\{ \frac{\partial^{2}l}{\partial x^{2}}\left(x_{0},y_{0}\right)\right\} -\cdots\nonumber \\
 &  & \qquad+(-1)^{n}\frac{(\Delta x)^{n}}{n!}\left\{ \frac{\partial^{n}l}{\partial x^{n}}\left(x_{0},y_{0}\right)\right\} +\cdots\label{eq:2-2}
\end{eqnarray}

and $\Delta x\rightarrow0$ in the equations (\ref{eq:1-3}) and (\ref{eq:1-4}),
then

\begin{eqnarray}
l\left(x_{0},y_{0}+\Delta y\right) & = & l\left(x_{0},y_{0}\right)+\Delta y\left\{ \frac{\partial l}{\partial y}\left(x_{0},y_{0}\right)\right\} +\frac{(\Delta y)^{2}}{2!}\left\{ \frac{\partial^{2}l}{\partial y^{2}}\left(x_{0},y_{0}\right)\right\} +\cdots\nonumber \\
 &  & \qquad+\frac{(\Delta y)^{n}}{n!}\left\{ \frac{\partial^{n}l}{\partial y^{n}}\left(x_{0},y_{0}\right)\right\} +\cdots\label{eq:2-3}\\
l\left(x_{0},y_{0}-\Delta y\right) & = & l\left(x_{0},y_{0}\right)-\Delta y\left\{ \frac{\partial l}{\partial y}\left(x_{0},y_{0}\right)\right\} +\frac{(\Delta y)^{2}}{2!}\left\{ \frac{\partial^{2}l}{\partial y^{2}}\left(x_{0},y_{0}\right)\right\} -\cdots\nonumber \\
 &  & \qquad+(-1)^{n}\frac{(\Delta y)^{n}}{n!}\left\{ \frac{\partial^{n}l}{\partial y^{n}}\left(x_{0},y_{0}\right)\right\} +\cdots\label{eq:2-4}
\end{eqnarray}

Sequentially, ignoring the $n^{th}$ order terms from eq. (\ref{eq:2-1})
and eq. (\ref{eq:2-2}), and denoting these new equations as $g_{1}^{(n)}\left(x_{0}+\Delta x,y_{0}\right)$
and $g_{1}^{(n)}\left(x_{0}-\Delta x,y_{0}\right)$ for $n=2,3,...$,
we obtain below set of equations.

\begin{eqnarray}
g_{1}^{(1)}\left(x_{0}+\Delta x,y_{0}\right)-g_{1}^{(1)}\left(x_{0}-\Delta x,y_{0}\right) & = & 2\Delta x\left\{ \frac{\partial l}{\partial x}\left(x_{0},y_{0}\right)\right\} \label{eq:2-5}\\
g_{1}^{(3)}\left(x_{0}+\Delta x,y_{0}\right)-g_{1}^{(3)}\left(x_{0}-\Delta x,y_{0}\right) & = & 2\Delta x\left\{ \frac{\partial l}{\partial x}\left(x_{0},y_{0}\right)\right\} +\frac{2(\Delta x)^{3}}{3!}\left\{ \frac{\partial^{3}l}{\partial x^{3}}\left(x_{0},y_{0}\right)\right\} \nonumber \\
\vdots &  & \vdots\nonumber \\
g_{1}^{(2k-1)}\left(x_{0}+\Delta x,y_{0}\right)-g_{1}^{(2k-1)}\left(x_{0}-\Delta x,y_{0}\right) & = & 2\Delta x\left\{ \frac{\partial l}{\partial x}\left(x_{0},y_{0}\right)\right\} +\frac{2(\Delta x)^{3}}{3!}\left\{ \frac{\partial^{3}l}{\partial x^{3}}\left(x_{0},y_{0}\right)\right\} +\nonumber \\
 &  & \cdots+\frac{2(\Delta x)^{(2k-1)}}{(2k-1)!}\left\{ \frac{\partial^{(2k-1)}l}{\partial x^{(2k-1)}}\left(x_{0},y_{0}\right)\right\} \nonumber 
\end{eqnarray}
\\
Therefore, 

\begin{eqnarray*}
\sum_{j=0}^{k}\left[\left\{ g_{1}^{(2j+1)}\left(x_{0}+\Delta x,y_{0}\right)-g_{1}^{(2j+1)}\left(x_{0}-\Delta x,y_{0}\right)\right\} -\left\{ g_{1}^{(2j-1)}\left(x_{0}+\Delta x,y_{0}\right)-g_{1}^{(2j-1)}\left(x_{0}-\Delta x,y_{0}\right)\right\} \right] & =
\end{eqnarray*}

\begin{eqnarray}
 &  & \sum_{j=1}^{k}\frac{2(\Delta x)^{(2j+1)}}{(2j+1)!}\left\{ \frac{\partial^{(2j+1)}l}{\partial x^{(2j+1)}}\left(x_{0},y_{0}\right)\right\} \label{eq:2-6}\\
 & = & \left\{ g_{1}^{(2k+1)}\left(x_{0}+\Delta x,y_{0}\right)-g_{1}^{(2k+1)}\left(x_{0}-\Delta x,y_{0}\right)\right\} -\left\{ g_{1}^{(1)}\left(x_{0}+\Delta x,y_{0}\right)-g_{1}^{(1)}\left(x_{0}-\Delta x,y_{0}\right)\right\} \nonumber 
\end{eqnarray}

Sequentially, ignoring the $n^{th}$ order terms from eq. (\ref{eq:2-3})
and eq. (\ref{eq:2-4}), and denoting these new equations as $g_{2}^{(n)}\left(x_{0},y_{0}+\Delta y\right)$
and $g_{2}^{(n)}\left(x_{0},y_{0}+\Delta y\right)$ for $n=2,3,...$,
we obtain below set of equations.

\begin{eqnarray}
g_{2}^{(1)}\left(x_{0},y_{0}+\Delta y\right)-g_{2}^{(1)}\left(x_{0},y_{0}-\Delta y\right) & = & 2\Delta y\left\{ \frac{\partial l}{\partial y}\left(x_{0},y_{0}\right)\right\} \label{eq:2-7}\\
g_{2}^{(3)}\left(x_{0},y_{0}+\Delta y\right)-g_{2}^{(3)}\left(x_{0},y_{0}-\Delta y\right) & = & 2\Delta y\left\{ \frac{\partial l}{\partial y}\left(x_{0},y_{0}\right)\right\} +\frac{2(\Delta y)^{3}}{3!}\left\{ \frac{\partial^{3}l}{\partial y^{3}}\left(x_{0},y_{0}\right)\right\} \nonumber \\
\vdots &  & \vdots\nonumber \\
g_{2}^{(2k-1)}\left(x_{0},y_{0}+\Delta y\right)-g_{2}^{(2k-1)}\left(x_{0},y_{0}-\Delta y\right) & = & 2\Delta y\left\{ \frac{\partial l}{\partial y}\left(x_{0},y_{0}\right)\right\} +\frac{2(\Delta y)^{3}}{3!}\left\{ \frac{\partial^{3}l}{\partial y^{3}}\left(x_{0},y_{0}\right)\right\} +\nonumber \\
 &  & \cdots+\frac{2(\Delta y)^{(2k-1)}}{(2k-1)!}\left\{ \frac{\partial^{(2k-1)}l}{\partial y^{(2k-1)}}\left(x_{0},y_{0}\right)\right\} \nonumber 
\end{eqnarray}

Therefore, 

\begin{eqnarray*}
\sum_{j=0}^{k}\left[\left\{ g_{2}^{(2j+1)}\left(x_{0},y_{0}+\Delta y\right)-g_{2}^{(2j+1)}\left(x_{0},y_{0}+\Delta y\right)\right\} -\left\{ g_{2}^{(2j-1)}\left(x_{0},y_{0}+\Delta y\right)-g_{2}^{(2j-1)}\left(x_{0},y_{0}+\Delta y\right)\right\} \right] & =
\end{eqnarray*}

\begin{eqnarray}
 &  & \sum_{j=1}^{k}\frac{2(\Delta y)^{(2j+1)}}{(2j+1)!}\left\{ \frac{\partial^{(2j+1)}l}{\partial y^{(2j+1)}}\left(x_{0},y_{0}\right)\right\} \label{eq:2-8}\\
 & = & \left\{ l^{(2k+1)}\left(x_{0},y_{0}+\Delta y\right)-l^{(2k+1)}\left(x_{0},y_{0}+\Delta y\right)\right\} -\left\{ l^{(1)}\left(x_{0},y_{0}+\Delta y\right)-l^{(1)}\left(x_{0},y_{0}-\Delta y\right)\right\} \nonumber 
\end{eqnarray}

Now substituting the eqs. (\ref{eq:3.8,1 for mu-3}) and (\ref{eq:3.10,do for mu-3})
in the eq. (\ref{eq:mu-3}), we get

\begin{eqnarray}
\mu\left(x_{0}+\Delta x,y_{0}\right) & = & -\frac{1}{\left[l\left(x_{0},y_{0}\right)+\sum_{j=1}^{\infty}\frac{(\Delta x)^{j}}{j!}\left\{ \frac{\partial^{j}l}{\partial x^{j}}\left(x_{0},y_{0}\right)\right\} \right]}\sum_{j=0}^{\infty}\frac{(\Delta x)^{j}}{j!}\left\{ \frac{\partial^{j+1}l}{\partial x^{j+1}}\left(x_{0},y_{0}\right)\right\} \label{eq:mu-(x+deltax0)}
\end{eqnarray}

and, substituting the eqs. (\ref{eq:3.9,1 for mu-4}) and (\ref{eq:3.11,do for mu-4})
in the eq. (\ref{eq:mu-4}), we get

\begin{eqnarray}
\mu\left(x_{0},y_{0}+\Delta y\right) & = & -\frac{1}{\left[l\left(x_{0},y_{0}\right)+\sum_{j=1}^{\infty}\frac{(\Delta y)^{j}}{j!}\left\{ \frac{\partial^{j}l}{\partial y^{j}}\left(x_{0},y_{0}\right)\right\} \right]}\sum_{j=0}^{\infty}\frac{(\Delta y)^{j}}{j!}\left\{ \frac{\partial^{j+1}l}{\partial y^{j+1}}\left(x_{0},y_{0}\right)\right\} \label{eq:mu(x0,y+deltay)}
\end{eqnarray}

Using eq. (\ref{eq:2-5}) and eq. (\ref{eq:2-7}), we obtain forces
of mortalities for two variables as follows:

\begin{eqnarray*}
\mu\left(x_{0,}y_{0}\right) & = & -\frac{\left\{ g_{1}^{(1)}\left(x_{0}+\Delta x,y_{0}\right)-g_{1}^{(1)}\left(x_{0}-\Delta x,y_{0}\right)\right\} }{2\Delta xl\left(x_{0},y_{0}\right)}\\
\mu\left(x_{0,}y_{0}\right) &  & -\frac{\left\{ g_{2}^{(1)}\left(x_{0},y_{0}+\Delta y\right)-g_{2}^{(1)}\left(x_{0},y_{0}-\Delta y\right)\right\} }{2\Delta yl\left(x_{0},y_{0}\right)}
\end{eqnarray*}

Hence the derivative forces of mortalities are as follows:

\begin{eqnarray}
\frac{\partial\mu\left(x_{,}y\right)}{\partial x} & = & \lim_{\Delta x\rightarrow0}\frac{1}{\Delta x}\left[\frac{\left\{ g_{1}^{(1)}\left(x_{0}+\Delta x,y_{0}\right)-g_{1}^{(1)}\left(x_{0}-\Delta x,y_{0}\right)\right\} }{2\Delta xl\left(x_{0},y_{0}\right)}-\frac{\sum_{j=0}^{\infty}\frac{(\Delta x)^{j}}{j!}\left\{ \frac{\partial^{j+1}l}{\partial x^{j+1}}\left(x_{0},y_{0}\right)\right\} }{\left[l\left(x_{0},y_{0}\right)+\sum_{j=1}^{\infty}\frac{(\Delta x)^{j}}{j!}\frac{\partial^{j}l}{\partial x^{j}}\left(x_{0},y_{0}\right)\right]}\right]\nonumber \\
\label{eq:mu-1-1}\\
\frac{\partial\mu\left(x_{,}y\right)}{\partial y} & = & \lim_{\Delta y\rightarrow0}\frac{1}{\Delta y}\left[\frac{\left\{ g_{2}^{(1)}\left(x_{0},y_{0}+\Delta y\right)-g_{2}^{(1)}\left(x_{0},y_{0}-\Delta y\right)\right\} }{2\Delta yl\left(x_{0},y_{0}\right)}-\frac{\sum_{j=0}^{\infty}\frac{(\Delta y)^{j}}{j!}\left\{ \frac{\partial^{j+1}l}{\partial y^{j+1}}\left(x_{0},y_{0}\right)\right\} }{\left[l\left(x_{0},y_{0}\right)+\sum_{j=1}^{\infty}\frac{(\Delta y)^{j}}{j!}\left\{ \frac{\partial^{j}l}{\partial y^{j}}\left(x_{0},y_{0}\right)\right\} \right]}\right]\nonumber \\
\label{eq:mu-2-1}
\end{eqnarray}

\section{Analysis of third order equations}

Suppose $s$ be the function number of survivors at age $x$ with
two more influencing variables $y$ and $z,$ then the three variable
survival function $s\left(x,y,z\right)$ evaluated at $\left\{ x_{o},y_{0},z_{0}\right\} $
can be written as 

\begin{eqnarray}
s\left(x_{0}+\Delta x,y_{0}+\Delta y,z_{0}+\Delta z\right) & = & s\left(x_{0},y_{0},z_{0}\right)+\Delta x\left\{ \frac{\partial s}{\partial x}\left(x_{0},y_{0},z_{0}\right)\right\} +\Delta y\left\{ \frac{\partial s}{\partial y}\left(x_{0},y_{0},z_{0}\right)\right\} \nonumber \\
 &  & +\Delta z\left\{ \frac{\partial s}{\partial z}\left(x_{0},y_{0},z_{0}\right)\right\} +\frac{\Delta x^{2}}{2}\left\{ \frac{\partial^{2}s}{\partial x^{2}}\left(x_{0},y_{0},z_{0}\right)\right\} \nonumber \\
 &  & +\frac{\Delta y^{2}}{2}\left\{ \frac{\partial^{2}s}{\partial y^{2}}\left(x_{0},y_{0},z_{0}\right)\right\} +\frac{\Delta z^{2}}{2}\left\{ \frac{\partial^{2}s}{\partial z^{2}}\left(x_{0},y_{0},z_{0}\right)\right\} \nonumber \\
 &  & +\Delta x\Delta y\left\{ \frac{\partial^{2}s}{\partial x\partial y}\left(x_{0},y_{0},z_{0}\right)\right\} +\Delta x\Delta z\left\{ \frac{\partial^{2}s}{\partial x\partial z}\left(x_{0},y_{0},z_{0}\right)\right\} \nonumber \\
 &  & +\Delta y\Delta z\left\{ \frac{\partial^{2}s}{\partial y\partial z}\left(x_{0},y_{0},z_{0}\right)\right\} +\nonumber \\
 &  & \cdots+\sum_{n=k}^{\infty}\left[\sum_{n_{1},n_{2},n_{3}}\frac{1}{n_{1}!n_{2}!n_{3}!}\frac{\partial^{k}s}{\partial x^{n_{1}}\partial y^{n_{2}}\partial z^{n_{3}}}\Delta x^{n_{1}}\Delta y^{n_{2}}\Delta z^{n_{3}}\right]\label{eq:Threevariables-1}
\end{eqnarray}

Similarly, 

\begin{eqnarray}
s\left(x_{0}-\Delta x,y_{0}-\Delta y,z_{0}-\Delta z\right) & = & s\left(x_{0},y_{0},z_{0}\right)-\Delta x\left\{ \frac{\partial s}{\partial x}\left(x_{0},y_{0},z_{0}\right)\right\} -\Delta y\left\{ \frac{\partial s}{\partial y}\left(x_{0},y_{0},z_{0}\right)\right\} \nonumber \\
 &  & -\Delta z\left\{ \frac{\partial s}{\partial z}\left(x_{0},y_{0},z_{0}\right)\right\} +\frac{\Delta x^{2}}{2}\left\{ \frac{\partial^{2}s}{\partial x^{2}}\left(x_{0},y_{0},z_{0}\right)\right\} \nonumber \\
 &  & +\frac{\Delta y^{2}}{2}\left\{ \frac{\partial^{2}s}{\partial y^{2}}\left(x_{0},y_{0},z_{0}\right)\right\} +\frac{\Delta z^{2}}{2}\left\{ \frac{\partial^{2}s}{\partial z^{2}}\left(x_{0},y_{0},z_{0}\right)\right\} \nonumber \\
 &  & +\Delta x\Delta y\left\{ \frac{\partial^{2}s}{\partial x\partial y}\left(x_{0},y_{0},z_{0}\right)\right\} +\Delta x\Delta z\left\{ \frac{\partial^{2}s}{\partial x\partial z}\left(x_{0},y_{0},z_{0}\right)\right\} \nonumber \\
 &  & +\Delta y\Delta z\left\{ \frac{\partial^{2}s}{\partial y\partial z}\left(x_{0},y_{0},z_{0}\right)\right\} +\nonumber \\
 &  & \cdots+(-1)^{k}\sum_{n=k}^{\infty}\left[\sum_{n_{1},n_{2},n_{3}}\frac{1}{n_{1}!n_{2}!n_{3}!}\frac{\partial^{k}s}{\partial x^{n_{1}}\partial y^{n_{2}}\partial z^{n_{3}}}\Delta x^{n_{1}}\Delta y^{n_{2}}\Delta z^{n_{3}}\right]\label{eq:threevariables-2}
\end{eqnarray}

Here $k=n_{1}+n_{2}+n_{3}.$ Let $h^{(n)}\left(x_{0}+\Delta x,y_{0}+\Delta y,z_{0}+\Delta z\right)$
and $h^{(n)}\left(x_{0}-\Delta x,y_{0}-\Delta y,z_{0}-\Delta z\right)$
denote equations after ignoring the terms with derivatives beginning
from the $(n+1)^{th}$ order ($n=1,2,...)$ in the eqs. (\ref{eq:Threevariables-1})
and (\ref{eq:threevariables-2}). We will obtain following equations:

\begin{eqnarray*}
\left.\begin{array}{c}
h^{(1)}\left(x_{0}+\Delta x,y_{0}+\Delta y,z_{0}+\Delta z\right)\\
-h^{(1)}\left(x_{0}-\Delta x,y_{0}-\Delta y,z_{0}-\Delta z\right)
\end{array}\right\}  & = & 2\Delta x\left\{ \frac{\partial s}{\partial x}\left(x_{0},y_{0},z_{0}\right)\right\} +2\Delta y\left\{ \frac{\partial s}{\partial y}\left(x_{0},y_{0},z_{0}\right)\right\} \\
 &  & +2\Delta z\left\{ \frac{\partial s}{\partial z}\left(x_{0},y_{0},z_{0}\right)\right\} \\
\vdots & \vdots & \vdots
\end{eqnarray*}

\begin{eqnarray*}
\left.\begin{array}{c}
h^{(2k-1)}\left(x_{0}+\Delta x,y_{0}+\Delta y,z_{0}+\Delta z\right)\\
-h^{(2k-1)}\left(x_{0}-\Delta x,y_{0}-\Delta y,z_{0}-\Delta z\right)
\end{array}\right\}  & = & 2\Delta x\left\{ \frac{\partial s}{\partial x}\left(x_{0},y_{0},z_{0}\right)\right\} +2\Delta y\left\{ \frac{\partial s}{\partial y}\left(x_{0},y_{0},z_{0}\right)\right\} \\
 &  & +2\Delta z\left\{ \frac{\partial s}{\partial z}\left(x_{0},y_{0},z_{0}\right)\right\} +\\
 &  & \cdots+\sum_{n_{1},n_{2},n_{3}}\frac{2}{n_{1}!n_{2}!n_{3}!}\frac{\partial^{(2k-1)}s}{\partial x^{n_{1}}\partial y^{n_{2}}\partial z^{n_{3}}}\Delta x^{n_{1}}\Delta y^{n_{2}}\Delta z^{n_{3}}\\
 &  & \left(\mbox{here }2k-1=n_{1}+n_{2}+n_{3}\right)\\
\vdots & \vdots & \vdots
\end{eqnarray*}

Therefore,

\begin{eqnarray*}
\sum_{j=1}^{k}\left[\left\{ h^{(2k+1)}\left(x_{0}+\Delta x,y_{0}+\Delta y,z_{0}+\Delta z\right)-h^{(2k+1)}\left(x_{0}-\Delta x,y_{0}-\Delta y,z_{0}-\Delta z\right)\right\} \right.\\
-\left.\left\{ g^{(2k-1)}\left(x_{0}+\Delta x,y_{0}+\Delta y,z_{0}+\Delta z\right)-g^{(2k-1)}\left(x_{0}-\Delta x,y_{0}-\Delta y,z_{0}+\Delta z\right)\right\} \right] & =
\end{eqnarray*}

\begin{eqnarray*}
 &  & \sum_{n_{1},n_{2},n_{3}}\frac{2}{n_{1}!n_{2}!n_{3}!}\frac{\partial^{3}s}{\partial x^{n_{1}}\partial y^{n_{2}}\partial z^{n_{3}}}\Delta x^{n_{1}}\Delta y^{n_{2}}\Delta z^{n_{3}}\\
 &  & \left(\mbox{here }3=n_{1}+n_{2}+n_{3}\right)
\end{eqnarray*}

\begin{eqnarray*}
 &  & =+\cdots+\sum_{n_{1},n_{2},n_{3}}\frac{2}{n_{1}!n_{2}!n_{3}!}\frac{\partial^{(2k+1)}s}{\partial x^{n_{1}}\partial y^{n_{2}}\partial z^{n_{3}}}\Delta x^{n_{1}}\Delta y^{n_{2}}\Delta z^{n_{3}}\\
 &  & \left(\mbox{here }2k+1=n_{1}+n_{2}+n_{3}\right)
\end{eqnarray*}

We define following three force of mortality functions evaluated at
$\left(x_{0},y_{0},z_{0}\right)$:

\begin{eqnarray*}
\mu_{x}\left(x_{0,}y_{0},z_{0}\right) & = & -\frac{1}{s\left(x_{0},y_{0},z_{0}\right)}\frac{\partial s\left(x_{0},y_{0},z_{0}\right)}{\partial x}\\
\mu_{y}\left(x_{0,}y_{0},z_{0}\right) & = & -\frac{1}{s\left(x_{0},y_{0},z_{0}\right)}\frac{\partial s\left(x_{0},y_{0},z_{0}\right)}{\partial y}\\
\mu_{z}\left(x_{0,}y_{0},z_{0}\right) & = & -\frac{1}{s\left(x_{0},y_{0},z_{0}\right)}\frac{\partial s\left(x_{0},y_{0},z_{0}\right)}{\partial z}
\end{eqnarray*}

and further we define three functions of forces of mortality as follows:

\begin{eqnarray*}
\mu\left(x_{0}+\Delta x,y_{0},z_{0}\right) & = & -\frac{1}{s\left(x_{0}+\Delta x,y_{0},z_{0}\right)}\frac{\partial s\left(x_{0}+\Delta x,y_{0},z_{0}\right)}{\partial x}\\
\mu\left(x_{0},y_{0}+\Delta y,z_{0}\right) & = & -\frac{1}{s\left(x_{0},y_{0}+\Delta y,z_{0}\right)}\frac{\partial s\left(x_{0},y_{0}+\Delta y,z_{0}\right)}{\partial y}\\
\mu\left(x_{0},y_{0},z_{0}+\Delta z\right) & = & -\frac{1}{s\left(x_{0},y_{0},z_{0}+\Delta z\right)}\frac{\partial s\left(x_{0},y_{0},z_{0}+\Delta z\right)}{\partial z}
\end{eqnarray*}

Using the above definitions, we obtain following rates evaluated at
$\left(x_{0},y_{0},z_{0}\right)$:

\begin{eqnarray*}
\frac{\partial\mu\left(x,y,z\right)}{\partial x} & = & \lim_{\Delta x\rightarrow0}\frac{\mu\left(x_{0}+\Delta x,y_{0},z_{0}\right)-\mu\left(x_{0},y_{0},z_{0}\right)}{\Delta x}\\
\frac{\partial\mu\left(x,y,z\right)}{\partial y} & = & \lim_{\Delta y\rightarrow0}\frac{\mu\left(x_{0},y_{0}+\Delta y,z_{0}\right)-\mu\left(x_{0},y_{0},z_{0}\right)}{\Delta y}\\
\frac{\partial\mu\left(x,y,z\right)}{\partial z} & = & \lim_{\Delta z\rightarrow0}\frac{\mu\left(x_{0},y_{0},z_{0}+\Delta z\right)-\mu\left(x_{0},y_{0},z_{0}\right)}{\Delta z}
\end{eqnarray*}

By taking pairs of limits$(\Delta y\rightarrow0$, $\Delta z\rightarrow0)$,
$(\Delta x\rightarrow0$, $\Delta z\rightarrow0)$, and $(\Delta x\rightarrow0$,
$\Delta y\rightarrow0)$, separately in the equation (\ref{eq:Threevariables-1}),
we obtain following three equations:

\begin{eqnarray}
s\left(x_{0}+\Delta x,y_{0},z_{0}\right) & = & s\left(x_{0},y_{0},z_{0}\right)+\Delta x\left\{ \frac{\partial s}{\partial x}\left(x_{0},y_{0},z_{0}\right)\right\} +\cdots+\frac{\Delta x^{n}}{n!}\left\{ \frac{\partial^{n}s}{\partial x^{n}}\left(x_{0},y_{0},z_{0}\right)\right\} +\nonumber \\
\label{eq:s1}
\end{eqnarray}

\begin{eqnarray}
s\left(x_{0},y_{0}+\Delta y,z_{0}\right) & = & s\left(x_{0},y_{0},z_{0}\right)+\Delta y\left\{ \frac{\partial s}{\partial y}\left(x_{0},y_{0},z_{0}\right)\right\} +\cdots+\frac{\Delta y^{n}}{n!}\left\{ \frac{\partial^{n}s}{\partial y^{n}}\left(x_{0},y_{0},z_{0}\right)\right\} +\cdots\nonumber \\
\label{eq:s2}
\end{eqnarray}

\begin{eqnarray}
s\left(x_{0},y_{0},z_{0}+\Delta z\right) & = & s\left(x_{0},y_{0},z_{0}\right)+\Delta z\left\{ \frac{\partial s}{\partial z}\left(x_{0},y_{0},z_{0}\right)\right\} +\cdots+\frac{\Delta z^{n}}{n!}\left\{ \frac{\partial^{n}s}{\partial z^{n}}\left(x_{0},y_{0},z_{0}\right)\right\} +\cdots\nonumber \\
\label{eq:s3}
\end{eqnarray}

By taking pairs of limits$(\Delta y\rightarrow0$, $\Delta z\rightarrow0)$,
$(\Delta x\rightarrow0$, $\Delta z\rightarrow0)$, and $(\Delta x\rightarrow0$,
$\Delta y\rightarrow0)$, separately in the equation (\ref{eq:threevariables-2}),
we obtain following three equations:

\begin{eqnarray}
s\left(x_{0}-\Delta x,y_{0},z_{0}\right) & = & s\left(x_{0},y_{0},z_{0}\right)-\Delta x\left\{ \frac{\partial s}{\partial x}\left(x_{0},y_{0},z_{0}\right)\right\} +\cdots+(-1)^{n}\frac{\Delta x^{n}}{n!}\left\{ \frac{\partial^{n}s}{\partial x^{n}}\left(x_{0},y_{0},z_{0}\right)\right\} +\cdots\nonumber \\
\end{eqnarray}

\begin{eqnarray}
s\left(x_{0},y_{0}-\Delta y,z_{0}\right) & = & s\left(x_{0},y_{0},z_{0}\right)-\Delta y\left\{ \frac{\partial s}{\partial y}\left(x_{0},y_{0},z_{0}\right)\right\} +\cdots+(-1)^{n}\frac{\Delta y^{n}}{n!}\left\{ \frac{\partial^{n}s}{\partial y^{n}}\left(x_{0},y_{0},z_{0}\right)\right\} +\cdots\nonumber \\
\end{eqnarray}

\begin{eqnarray}
s\left(x_{0},y_{0},z_{0}-\Delta z\right) & = & s\left(x_{0},y_{0},z_{0}\right)-\Delta z\left\{ \frac{\partial s}{\partial z}\left(x_{0},y_{0},z_{0}\right)\right\} +\cdots+(-1)^{n}\frac{\Delta z^{n}}{n!}\left\{ \frac{\partial^{n}s}{\partial z^{n}}\left(x_{0},y_{0},z_{0}\right)\right\} +\cdots\nonumber \\
\end{eqnarray}

Suppose $h_{1}^{n}\left(x_{0}+\Delta x,y_{0},z_{0}\right)$, $h_{1}^{n}\left(x_{0}-\Delta x,y_{0},z_{0}\right)$;
$h_{2}^{n}\left(x_{0},y_{0}+\Delta y,z_{0}\right)$, $h_{2}^{n}\left(x_{0},y_{0}-\Delta y,z_{0}\right)$
and

$h_{3}^{n}\left(x_{0},y_{0},z_{0}+\Delta z\right)$ , $h_{3}^{n}\left(x_{0},y_{0},z_{0}-\Delta z\right)$
for $n=1,2,3,...$ denote the functions by ignoring the terms from
the order $(n+1)$ in the pairs of equations (\ref{eq:s1}), (\ref{eq:s2});
(\ref{eq:s2}), (\ref{eq:s3}) and (\ref{eq:s3}), (\ref{eq:s3}),
then we will obtain following three series of sequences of difference
functions:

\begin{eqnarray*}
h_{1}^{(1)}\left(x_{0}+\Delta x,y_{0},z_{0}\right)-h_{1}^{(1)}\left(x_{0}-\Delta x,y_{0},z_{0}\right) & = & 2\Delta x\left\{ \frac{\partial s}{\partial x}\left(x_{0},y_{0},z_{0}\right)\right\} \\
h_{1}^{(3)}\left(x_{0}+\Delta x,y_{0},z_{0}\right)-h_{1}^{(3)}\left(x_{0}-\Delta x,y_{0},z_{0}\right) & = & 2\Delta x\left\{ \frac{\partial s}{\partial x}\left(x_{0},y_{0},z_{0}\right)\right\} +\frac{2\Delta x^{3}}{3!}\left\{ \frac{\partial^{3}s}{\partial x^{3}}\left(x_{0},y_{0},z_{0}\right)\right\} 
\end{eqnarray*}

\begin{eqnarray*}
\vdots & \vdots & \vdots\\
h_{1}^{(2k-1)}\left(x_{0}+\Delta x,y_{0},z_{0}\right)-h_{1}^{(2k-1)}\left(x_{0}-\Delta x,y_{0},z_{0}\right) & = & 2\Delta x\left\{ \frac{\partial s}{\partial x}\left(x_{0},y_{0},z_{0}\right)\right\} +\frac{2\Delta x^{3}}{3!}\left\{ \frac{\partial^{3}s}{\partial x^{3}}\left(x_{0},y_{0},z_{0}\right)\right\} +\\
 &  & \cdots+\frac{2\Delta x^{(2k-1)}}{(2k-1)!}\left\{ \frac{\partial^{(2k-1)}s}{\partial x^{(2k-1)}}\left(x_{0},y_{0},z_{0}\right)\right\} \\
\vdots & \vdots & \vdots
\end{eqnarray*}

Therefore, 

\begin{eqnarray}
\sum_{j=1}^{k}\left[\begin{array}{c}
\left\{ h_{1}^{(2j+1)}\left(x_{0}+\Delta x,y_{0},z_{0}\right)-h_{1}^{(2j+1)}\left(x_{0}-\Delta x,y_{0},z_{0}\right)\right\} -\\
\left\{ h_{1}^{(2j-1)}\left(x_{0}+\Delta x,y_{0},z_{0}\right)-h_{1}^{(2j-1)}\left(x_{0}-\Delta x,y_{0},z_{0}\right)\right\} 
\end{array}\right] & = & \frac{2\Delta x^{3}}{3!}\left\{ \frac{\partial^{3}s}{\partial x^{3}}\left(x_{0},y_{0},z_{0}\right)\right\} +\nonumber \\
 &  & \cdots+\frac{2\Delta x^{(2k+1)}}{(2k+1)!}\left\{ \frac{\partial^{(2k+1)}s}{\partial x^{(2k+1)}}\left(x_{0},y_{0},z_{0}\right)\right\} \label{eq:series-1}
\end{eqnarray}

\begin{eqnarray*}
h_{2}^{(1)}\left(x_{0},y_{0}+\Delta y,z_{0}\right)-h_{2}^{(1)}\left(x_{0},y_{0}-\Delta y,z_{0}\right) & = & 2\Delta y\left\{ \frac{\partial s}{\partial y}\left(x_{0},y_{0},z_{0}\right)\right\} \\
h_{2}^{(3)}\left(x_{0},y_{0}+\Delta y,z_{0}\right)-h_{2}^{(3)}\left(x_{0},y_{0}-\Delta y,z_{0}\right) & = & 2\Delta y\left\{ \frac{\partial s}{\partial y}\left(x_{0},y_{0},z_{0}\right)\right\} +\frac{2\Delta y^{3}}{3!}\left\{ \frac{\partial^{3}s}{\partial y^{3}}\left(x_{0},y_{0},z_{0}\right)\right\} 
\end{eqnarray*}

\begin{eqnarray*}
\vdots & \vdots & \vdots\\
h_{2}^{(2k-1)}\left(x_{0},y_{0}+\Delta y,z_{0}\right)-h_{2}^{(2k-1)}\left(x_{0},y_{0}-\Delta y,z_{0}\right) & = & 2\Delta y\left\{ \frac{\partial s}{\partial y}\left(x_{0},y_{0},z_{0}\right)\right\} ++\frac{2\Delta y^{3}}{3!}\left\{ \frac{\partial^{3}s}{\partial y^{3}}\left(x_{0},y_{0},z_{0}\right)\right\} \\
 &  & \cdots+\frac{2\Delta y^{2k-1}}{(2k-1)!}\left\{ \frac{\partial^{2k-1}s}{\partial y^{2k-1}}\left(x_{0},y_{0},z_{0}\right)\right\} \\
\vdots & \vdots & \vdots
\end{eqnarray*}

Therefore, 

\begin{eqnarray}
\sum_{j=1}^{k}\left[\begin{array}{c}
\left\{ h_{2}^{(2j+1)}\left(x_{0},y_{0}+\Delta y,z_{0}\right)-h_{2}^{(2j+1)}\left(x_{0},y_{0}-\Delta y,z_{0}\right)\right\} -\\
\left\{ h_{2}^{(2j-1)}\left(x_{0},y_{0}+\Delta y,z_{0}\right)-h_{2}^{(2j-1)}\left(x_{0},y_{0}-\Delta y,z_{0}\right)\right\} 
\end{array}\right] & = & \frac{2\Delta y^{3}}{3!}\left\{ \frac{\partial^{3}s}{\partial y^{3}}\left(x_{0},y_{0},z_{0}\right)\right\} +\nonumber \\
 &  & \cdots+\frac{2\Delta y^{2k+1}}{(2k+1)!}\left\{ \frac{\partial^{2k+1}s}{\partial y^{2k+1}}\left(x_{0},y_{0},z_{0}\right)\right\} \label{eq:series-2}
\end{eqnarray}

\begin{eqnarray*}
h_{3}^{(1)}\left(x_{0},y_{0},z_{0}+\Delta z\right)-h_{3}^{(1)}\left(x_{0},y_{0},z_{0}-\Delta z\right) & = & 2\Delta z\left\{ \frac{\partial s}{\partial z}\left(x_{0},y_{0},z_{0}\right)\right\} \\
h_{3}^{(3)}\left(x_{0},y_{0},z_{0}+\Delta z\right)-h_{3}^{(3)}\left(x_{0},y_{0},z_{0}-\Delta z\right) & = & 2\Delta z\left\{ \frac{\partial s}{\partial z}\left(x_{0},y_{0},z_{0}\right)\right\} +\frac{2\Delta z^{3}}{3!}\left\{ \frac{\partial^{3}s}{\partial z^{3}}\left(x_{0},y_{0},z_{0}\right)\right\} 
\end{eqnarray*}

\begin{eqnarray*}
\vdots & \vdots & \vdots\\
h_{3}^{(2k-1)}\left(x_{0},y_{0},z_{0}+\Delta z\right)-h_{3}^{(2k-1)}\left(x_{0},y_{0},z_{0}-\Delta z\right) & = & 2\Delta z\left\{ \frac{\partial s}{\partial z}\left(x_{0},y_{0},z_{0}\right)\right\} +\frac{2\Delta z^{3}}{3!}\left\{ \frac{\partial^{3}s}{\partial z^{3}}\left(x_{0},y_{0},z_{0}\right)\right\} \\
 &  & \cdots+\frac{2\Delta z^{(2k-1)}}{(2k-1)!}\left\{ \frac{\partial^{(2k-1)}s}{\partial z^{(2k-1)}}\left(x_{0},y_{0},z_{0}\right)\right\} \\
\vdots & \vdots & \vdots
\end{eqnarray*}

Therefore,

\begin{eqnarray}
\sum_{j=1}^{k}\left[\begin{array}{c}
\left\{ h_{3}^{(2j+1)}\left(x_{0},y_{0},z_{0}+\Delta z\right)-h_{3}^{(2j+1)}\left(x_{0},y_{0},z_{0}-\Delta z\right)\right\} -\\
\left\{ h_{3}^{(2j-1)}\left(x_{0},y_{0},z_{0}+\Delta z\right)-h_{3}^{(2j-1)}\left(x_{0},y_{0},z_{0}-\Delta z\right)\right\} 
\end{array}\right] & = & \frac{2\Delta z^{3}}{3!}\left\{ \frac{\partial^{3}s}{\partial z^{3}}\left(x_{0},y_{0},z_{0}\right)\right\} +\nonumber \\
 &  & \cdots+\frac{2\Delta z^{(2k+1)}}{(2k+1)!}\left\{ \frac{\partial^{(2k+1)}s}{\partial z^{(2k+1)}}\left(x_{0},y_{0},z_{0}\right)\right\} \label{eq:series-3}
\end{eqnarray}

The rate of changes in the survival function with respect to one variable
and corresponding forces of mortalities for three variables can be
obtained using the following derivations.

\begin{eqnarray*}
\frac{\partial s\left(x_{0}+\Delta x,y_{0},z_{0}\right)}{\partial x} & = & \frac{\partial s\left(x_{0},y_{0},z_{0}\right)}{\partial x}+\Delta x\left\{ \frac{\partial^{2}s}{\partial x^{2}}\left(x_{0},y_{0},z_{0}\right)\right\} +\cdots+\frac{\Delta x^{n}}{n!}\left\{ \frac{\partial^{n+1}s}{\partial x^{n+1}}\left(x_{0},y_{0},z_{0}\right)\right\} +\cdots
\end{eqnarray*}

\begin{eqnarray*}
\frac{\partial s\left(x_{0},y_{0}+\Delta y,z_{0}\right)}{\partial y} & = & \frac{\partial s\left(x_{0},y_{0},z_{0}\right)}{\partial y}+\Delta y\left\{ \frac{\partial^{2}s}{\partial y^{2}}\left(x_{0},y_{0},z_{0}\right)\right\} +\cdots+\frac{\Delta y^{n}}{n!}\left\{ \frac{\partial^{n+1}s}{\partial y^{n+1}}\left(x_{0},y_{0},z_{0}\right)\right\} +\cdots
\end{eqnarray*}

\begin{eqnarray*}
\frac{\partial s\left(x_{0},y_{0},z_{0}+\Delta z\right)}{\partial z} & = & \frac{\partial s\left(x_{0},y_{0},z_{0}\right)}{\partial z}+\Delta z\left\{ \frac{\partial^{2}s}{\partial z^{2}}\left(x_{0},y_{0},z_{0}\right)\right\} +\cdots+\frac{\Delta z^{n}}{n!}\left\{ \frac{\partial^{n+1}s}{\partial z^{n+1}}\left(x_{0},y_{0},z_{0}\right)\right\} +\cdots
\end{eqnarray*}

\begin{eqnarray*}
\mu\left(x_{0}+\Delta x,y_{0},z_{0}\right) & = & -\frac{1}{s\left(x_{0},y_{0},z_{0}\right)+\sum_{i=1}^{\infty}\left[\frac{\Delta x^{i}}{i!}\left\{ \frac{\partial^{i}s}{\partial x^{i}}\left(x_{0},y_{0},z_{0}\right)\right\} \right]}\sum_{i=0}^{\infty}\left[\frac{\Delta x^{i}}{i!}\left\{ \frac{\partial^{i+1}s}{\partial x^{i+1}}\left(x_{0},y_{0},z_{0}\right)\right\} \right]
\end{eqnarray*}

\begin{eqnarray*}
\mu\left(x_{0},y_{0}+\Delta y,z_{0}\right) & = & -\frac{1}{s\left(x_{0},y_{0},z_{0}\right)+\sum_{i=1}^{\infty}\left[\frac{\Delta y^{i}}{i!}\left\{ \frac{\partial^{i}s}{\partial y^{i}}\left(x_{0},y_{0},z_{0}\right)\right\} \right]}\sum_{i=0}^{\infty}\left[\frac{\Delta y^{i}}{i!}\left\{ \frac{\partial^{i+1}s}{\partial y^{i+1}}\left(x_{0},y_{0},z_{0}\right)\right\} \right]
\end{eqnarray*}

\begin{eqnarray*}
\mu\left(x_{0},y_{0},z_{0}+\Delta z\right) & = & -\frac{1}{s\left(x_{0},y_{0},z_{0}\right)+\sum_{i=1}^{\infty}\left[\frac{\Delta z^{i}}{i!}\left\{ \frac{\partial^{i}s}{\partial z^{i}}\left(x_{0},y_{0},z_{0}\right)\right\} \right]}\sum_{i=0}^{\infty}\left[\frac{\Delta z^{i}}{i!}\left\{ \frac{\partial^{i+1}s}{\partial z^{i+1}}\left(x_{0},y_{0},z_{0}\right)\right\} \right]
\end{eqnarray*}

\begin{eqnarray*}
\frac{\partial\mu\left(x,y,z\right)}{\partial x} & = & \lim_{\Delta x\rightarrow0}\frac{1}{\Delta x}\left[\frac{1}{s\left(x_{0},y_{0},z_{0}\right)}\left\{ \frac{h_{1}^{(1)}\left(x_{0}+\Delta x,y_{0},z_{0}\right)-h_{1}^{(1)}\left(x_{0}-\Delta x,y_{0},z_{0}\right)}{2\Delta x}\right\} \right.\\
 &  & -\left.\left\{ \frac{\sum_{i=0}^{\infty}\left[\frac{\Delta x^{i}}{i!}\left\{ \frac{\partial^{i+1}s}{\partial x^{i+1}}\left(x_{0},y_{0},z_{0}\right)\right\} \right]}{s\left(x_{0},y_{0},z_{0}\right)+\sum_{i=1}^{\infty}\left[\frac{\Delta x^{i}}{i!}\left\{ \frac{\partial^{i}s}{\partial x^{i}}\left(x_{0},y_{0},z_{0}\right)\right\} \right]}\right\} \right]
\end{eqnarray*}

\begin{eqnarray*}
\frac{\partial\mu\left(x,y,z\right)}{\partial y} & = & \lim_{\Delta y\rightarrow0}\frac{1}{\Delta y}\left[\frac{1}{s\left(x_{0},y_{0},z_{0}\right)}\left\{ \frac{h_{y}^{(1)}\left(x_{0},y_{0}+\Delta y,z_{0}\right)-h_{2}^{(1)}\left(x_{0},y_{0}-\Delta y,z_{0}\right)}{2\Delta y}\right\} \right.\\
 &  & -\left.\left\{ \frac{\sum_{i=0}^{\infty}\left[\frac{\Delta y^{i}}{i!}\left\{ \frac{\partial^{i+1}s}{\partial y^{i+1}}\left(x_{0},y_{0},z_{0}\right)\right\} \right]}{s\left(x_{0},y_{0},z_{0}\right)+\sum_{i=1}^{\infty}\left[\frac{\Delta y^{i}}{i!}\left\{ \frac{\partial^{i}s}{\partial y^{i}}\left(x_{0},y_{0},z_{0}\right)\right\} \right]}\right\} \right]
\end{eqnarray*}

\begin{eqnarray*}
\frac{\partial\mu\left(x,y,z\right)}{\partial z} & = & \lim_{\Delta z\rightarrow0}\frac{1}{\Delta z}\left[\frac{1}{s\left(x_{0},y_{0},z_{0}\right)}\left\{ \frac{h_{3}^{(1)}\left(x_{0},y_{0},z_{0}+\Delta z\right)-h_{3}^{(1)}\left(x_{0},y_{0},z_{0}-\Delta z\right)}{2\Delta z}\right\} \right.\\
 &  & -\left.\left\{ \frac{\sum_{i=0}^{\infty}\left[\frac{\Delta y^{i}}{i!}\left\{ \frac{\partial^{i+1}s}{\partial y^{i+1}}\left(x_{0},y_{0},z_{0}\right)\right\} \right]}{s\left(x_{0},y_{0},z_{0}\right)+\sum_{i=1}^{\infty}\left[\frac{\Delta y^{i}}{i!}\left\{ \frac{\partial^{i}s}{\partial y^{i}}\left(x_{0},y_{0},z_{0}\right)\right\} \right]}\right\} \right]
\end{eqnarray*}

\section{Conclusions}

Our numerical examples and analytical derivations does encourage to
validate results obtained by univariate force of mortality with that
of bivariate and multivariate forces of mortality functions. Majority
of the mortality data analysed consider age as a predominant variable
\cite{Bebbington2007,Finkelstein2005,Garilov1991,TuljaLiBoeNature2000}
in forecasting and analysis. Some insect populations as well age is
considered as a predominant variable in mortality analysis \cite{Carey1992}.
Mortality data analysed does indicate that considering only one variable
in concluding the causes of decline could lead to incomplete hypothesis,
thus warrants further analysis.

\end{document}